\documentclass{amsart}
\usepackage{amssymb,amsbsy,amsmath,amsfonts,amssymb,amscd}
\usepackage{latexsym,eucal,exscale,epic,eepic,epsfig}
\usepackage[english,francais]{babel}
\usepackage{times}

\newcommand{\integers}{{\bf Z}}
\newcommand{\reals}{{\bf R}}
\newcommand{\complexes}{{\bf C}}
\newcommand{\eisenstein}{\EuScript{E}}

\newcommand{\chyper}{{\bf CH}}
\newcommand{\rproj}{{\bf RP}}
\newcommand{\rhyper}{{\bf RH}}

\newcommand{\cubicForms}{\EuScript{C}}
\newcommand{\moduli}{\EuScript{M}}

\newcommand{\calH}{\EuScript{H}}

\newcommand{\map}{\longrightarrow}
\newcommand{\comp}{\circ}

\newtheorem{theorem}{Theorem}[section]

\newtheorem{e-proposition}[theorem]{Proposition}

\newtheorem{e-definition}[theorem]{Definition\rm}

\begin{document}
\selectlanguage{english}

\title{Real Cubic Surfaces and Real Hyperbolic Geometry}
\author{Daniel Allcock}
\address{Department of Mathematics\\University of Texas\\Austin, TX 78712}
\email{allcock@math.utexas.edu}
\urladdr{http://www.math.utexas.edu/\textasciitilde allcock}
\author{James A. Carlson}
\address{Department of Mathematics\\University of Utah\\Salt
Lake City, UT 84112}
\email{carlson@math.utah.edu}
\urladdr{http://www.math.utah.edu/\textasciitilde carlson}
\author{Domingo Toledo}
\address{Department of Mathematics\\University of Utah\\Salt
Lake City, UT 84112}
\email{toledo@math.utah.edu}
\urladdr{http://www.math.utah.edu/\textasciitilde toledo}
\date{March 29, 2003}
\thanks{First author partly supported by NSF grant DMS~0070930.
Second and third authors partly supported by NSF grant DMS~9900543.}
%

%
\begin{abstract}{%
The moduli space of stable real cubic surfaces is the quotient of real
hyperbolic four-space by a discrete, nonarithmetic group.  The volume
of the moduli space is $37\pi^2/1080$ in the metric of constant
curvature $-1$.  Each of the five connected components of the moduli
space can be described as the quotient of real hyperbolic
four-space by a specific arithmetic group. We compute the volumes of these
components.  }\end{abstract}


\maketitle
\thispagestyle{empty}

\setcounter{section}{0}
\selectlanguage{english}

\section{Results}

In \cite{ACT:JGA} we showed that the moduli space of stable cubic
surfaces is the quotient of complex hyperbolic four-space by a certain
arithmetic group which we described explicitly.  The purpose of this
note is announce a corresponding result for real cubic surfaces:
the moduli space is a quotient of real hyperbolic four-space by an
explicit discrete group. The group, however, is not
arithmetic.  We also compute the volume of the moduli
space in its metric of curvature $-1$.  It is $37\pi^2/1080 =
(4\pi^2/3)(37/1440)$.  (The $4\pi^2/3$ is the ratio of the volume of
the unit 4-sphere to its Euler characteristic, which appears in the
Gauss-Bonnet theorem.)

By the moduli space $\moduli_0^\reals$ (resp. $\moduli_s^\reals$) we
mean the set $\cubicForms_0^\reals$ (resp. $\cubicForms_s^\reals$) of
cubic forms with real coefficients that define smooth (resp. stable)
surfaces, modulo the action of $GL(4,\reals)$.  By smooth we mean that
the set of complex points is smooth, and by stable we mean stable in
the sense of geometric invariant theory.  In this case, stable means
that the complex surface has no singularities besides nodes.  The
space $\moduli_0^\reals$ has five connected components (see
\cite{Segre}), which we denote by $\moduli_{0,j}^\reals$ for $j = 0,
1, \dots,4$.

For each component $\moduli_j^\reals$ of the moduli space we exhibit
an arithmetic lattice $P\Gamma_j \subset PO(4,1)$, a union $\Delta_j$
of two- and three-dimensional real hyperbolic subspaces of
$\rhyper^{\smash{4}}$, and an isomorphism
\begin{equation}
  \moduli_{0,j}^\reals \cong P\Gamma_j \backslash (\rhyper^4 - \Delta_j)
\label{mainequation}
\end{equation}
of real analytic orbifolds.  We give two concrete descriptions of the
$P\Gamma_j$, one arithmetic and one geometric.  First, $P\Gamma_j$ is
the projective orthogonal group of the integer quadratic form $ -x_0^2
+ m_1x_1^2 + \cdots + m_4x_4^2$, where $j$ of the $m_i$ are 3's and
the rest are 1's.  Second, $P\Gamma_j$ is, up to a group of order at
most two, the Coxeter group $W_j$ defined in Figure~1.  More
precisely, $P\Gamma_j$ is the semidirect product of $W_j$ by the group
of diagram automorphisms, which is either trivial or of order two.
Yoshida has treated the case $j=0$ in \cite{Yoshida}.  


\begin{figure}
\def\LLx{0}
\def\LLy{0}
\def\width{320}
\def\height{220}
\def\Wtwoheight{35}
\def\Wthreeheight{25}
\setlength{\unitlength}{1bp}
\leavevmode\hskip.7in
\begin{picture}(\width,\height)(\LLx,\LLy)
\put(15,\height){\makebox(0,0)[tl]{\input smallW0.tex
}}
\put(305,\height){\makebox(0,0)[tr]{\input smallW1.tex
}}
\put(0,\Wtwoheight){\makebox(0,0)[bl]{\input  smallW2.tex
}}
\put(\width,\Wthreeheight){\makebox(0,0)[br]{\input smallW3.tex
}}
\put(300,0){\makebox(0,0)[br]{\input smallW4.tex
}}
\end{picture}
\caption{Coxeter diagrams for the reflection subgroups $W_j$ of the
$P\Gamma_j$.  Each describes a polyhedron $C_j$ with one facet per
atom of the diagram.  The bonds indicate if/how pairs of facets meet:
an absent (resp. single, double, triple) bond represents an angle of
$\frac{\pi}{2}$ (resp $\frac{\pi}{3}$, $\frac{\pi}{4}$, $\frac{\pi}{6}$), and a dashed
(resp. heavy) bond represents ultraparallelism (resp. parallelism at
$\infty$).  $W_j$ is the group generated by reflections in the facets
of $C_j$, and $C_j$ is a fundamental domain for $W_j$.}
\end{figure}
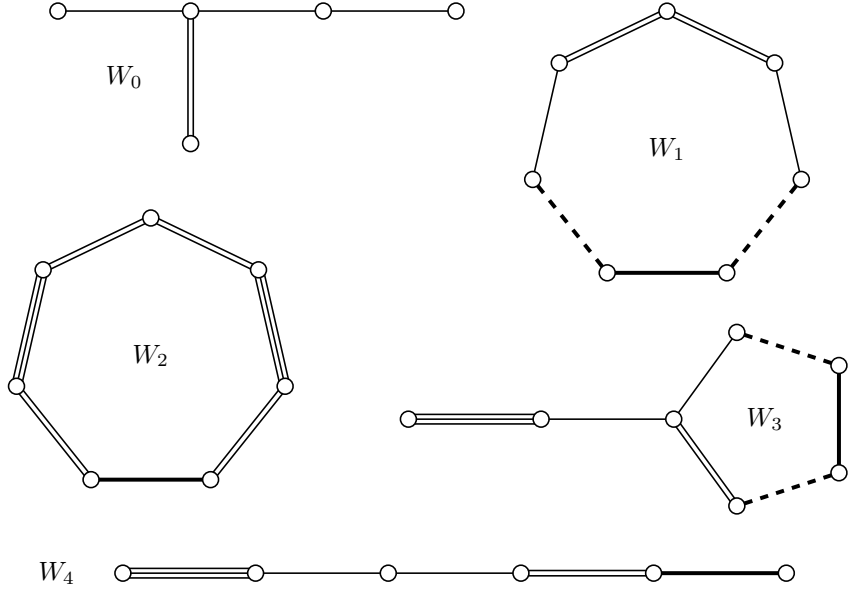


The points of $\Delta_j$ represent nodal surfaces which are limits of
smooth surfaces of type $j$.  Since a surface with a real node is a
limit of two different topological types of real surface, it is
natural to glue various pairs $\moduli_{0,j}^\reals$ and
$\moduli_{0,j'}^\reals$ together by identifying part of $\Delta_j$
with part of $\Delta_{j'}$.  Carrying this out in practice means
gluing certain faces of the polyhedra $C_j$ to each other and taking
care to deal with the diagram automorphisms.  A miracle
occurs and the result of these gluings turns out to be a quotient of
$\rhyper^{\smash{4}}$ in its own right:

\begin{theorem}
\label{maintheorem}
There is a nonarithmetic lattice 
$P\Gamma^\reals\subset PO(4,1)$, a union $\Delta$ of two- and
  three-dimensional hyperbolic subspaces of $\rhyper^4$, and
an isomorphism
\[
   \moduli_0^\reals \cong P\Gamma^\reals \backslash (\rhyper^4 - \Delta)
\]
of real analytic manifolds.  This identification extends to a
homeomorphism 
$$\moduli_s^\reals \cong P\Gamma^\reals \backslash \rhyper^4\;.$$
\end{theorem}

\smallskip
The clue to the nonarithmeticity is that $\moduli_s^\reals$ is
obtained by gluing together arithmetic orbifolds whose groups fall
into two commensurability classes.  In the spirit of
Gromov-Piatetski-Shapiro \cite{Gromov}, one expects the resulting
group to be non-arithmetic.  To prove the nonarithmeticity we use the
Galois-conjugation criterion in \cite{DM}.  Namely, it happens that
$P\Gamma^\reals$ preserves an integral quadratic form over
$\integers[\sqrt{3}]$ which has signature $(1,4)$ and whose
Galois conjugate has signature $(3,2)$.  We note that $P\Gamma^\reals$
is not a Coxeter group, even up to finite index, but it contains an
index two subgroup whose fundamental domain is a union of the $C_j$
and happens to be a Coxeter polyhedron.

The homeomorphism $\moduli_s^\reals\cong
P\Gamma^\reals\backslash\rhyper^4$ is not an orbifold isomorphism,
but it becomes one if the orbifold structure on
$P\Gamma^\reals\backslash\rhyper^4$ is suitably changed.  This can be
done explicitly enough to compute the orbifold fundamental group
$\pi_1^{orb}(\moduli_s^\reals)\cong\integers/2\times(\integers*\integers/2)$,
and to see that $\moduli_s^\reals$ is a bad orbifold in the sense of
Thurston.

The theory of Coxeter groups makes it easy to compute the orbifold
Euler characteristic of $W_j\backslash\rhyper^4$, and hence the
volume of this quotient.  Dividing by a factor of two if necessary, we
obtain the volume of $ P\Gamma_j\backslash\rhyper^4 $, which is the
volume of $\moduli_{0,j}^\reals$.  It follows that the hyperbolic
volume of $P\Gamma^\reals\backslash\rhyper^4$ is the sum of these
volumes.  The results are displayed in the table below.  For each $j$
we give the topology of that type of real cubic surface, the number of its
real lines, the orbifold fundamental group of $\moduli_{0,j}^\reals$,
and the orbifold Euler characteristic and volume of
$P\Gamma_j\backslash\rhyper^{\smash{4}}$.  $S_n$ and $D_\infty$ denote
symmetric and infinite dihedral groups.  Note that the component corresponding
to the simplest topology has the greatest volume, just over $40\%$ of
the total, and the component corresponding to surfaces with the most
real lines has the smallest volume.

\smallskip

\rlap{\begin{tabular}[t]{ccccrrr}
Type&Topology&\llap{Re}al Lin\rlap{es}&$\pi_1^{orb}(\moduli_{0,j}^\reals)$&\llap{Eu}ler char.&Volume& Fraction\\ 
\noalign{\vskip1.5pt}
\hline
\noalign{\vskip2pt}
0 & $\rproj^2 + 3\hbox{ handles}$& 27 &$S_5$& $1/1920$ & .00685 & 2.03\%\\
1 & $\rproj^2 + 2\hbox{ handles}$& 15  & $(S_3\times S_3) \rtimes \integers/2$& 1/288 &  .04569  &  13.51\%\\
2 & $\rproj^2 + 1\hbox{ handle}$ &7 & $(D_\infty\times D_\infty) \rtimes \integers/2$& 5/576 &      .11423 & 33.78\%\\
3 & $\rproj^2$ &3 &%
\newcommand{\braceheight}{12pt}
\newcommand{\braceloweramount}{6pt}
\newcommand{\diagloweramount}{9pt}
\smash{\lower\braceloweramount\hbox{$\left.\vrule width0pt height\braceheight\right\rbrace$}}%
\smash{\lower\diagloweramount\hbox{\begin{picture}(83,8)(296,396)
\setlength{\unitlength}{1bp}
\put(296,396){\makebox(0,0)[bl]{\epsfig{file=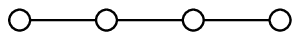, bbllx=296, bblly=396, bburx=379, bbury=404}}}
\put(312.5,400.0){\makebox(0,0)[b]{\kern 0pt\raise 3pt\hbox{$\infty$}}}
\end{picture}

}}%
& 1/96 &  .13708  & 40.54\%\\
4 & $\rproj^2 \cup S^2$ & 3 && 1/384 & .03427  & 10.14\%\\
\hline
\noalign{\vskip1pt}
&  &&  & 37/1440 &  .33813 & 100.00\% \\
\end{tabular}}

\section{About the proof}

The identification of the components of the moduli space with
quotients of real hyperbolic space depends on the construction of
\cite{ACT:CR}, \cite{ACT:JGA}.  Given a complex cubic surface $S$, let
$T$ be the triple cover of projective 3-space branched along $S$, and
let $(H^3(T), \sigma)$ denote the resulting special Hodge structure, where
$\sigma$ is the symmetry coming from the branched covering
transformation.  The period map which assigns to $S$ the class of
$(H^3(T), \sigma)$ defines an isomorphism between the moduli space of
stable cubic surfaces and $P\Gamma\backslash\chyper^4$.  Here
$P\Gamma$ is the projective automorphism group of the hermitian form
$h(x,y) = -x_0\bar y_0 + x_1\bar y_1 + \cdots + x_4\bar y_4$ on the
lattice $\Lambda = \eisenstein^{4,1}$, where $\eisenstein =
\integers[\sqrt[3]{1}]$.  The locus $\calH$ of $\chyper^4$
representing singular surfaces is the union of the orthogonal
complements of the norm~1 vectors of $\Lambda$.  In more detail, the
Hodge structure on $H^3(T)$, together with a choice of isomorphism
$i:H^3(T,\integers)\rightarrow\Lambda$ of Hermitian
$\eisenstein$-modules determines a complex line in $\Lambda_\complexes
= \Lambda\otimes_\eisenstein\complexes \cong \complexes^{4,1}$ which
is negative for $h$.  Thus $L$ is a point of $\chyper^4$, well defined
up to the action of $P\Gamma$.

We call an antilinear involution (``anti-involution'') of $\chyper^4$
integral if it arises from an anti-involution of $\Lambda$.  We write
$K_0$ for the set of all pairs $(L,\chi)$ where $L\in\chyper^4-\calH$
and $\chi$ is an integral anti-involution that preserves $L$.  
If the surface $S$ is defined by an equation with real coefficients,
then complex conjugation $\kappa(X_0,\ldots,X_4)
= (\bar X_0,\ldots, \bar X_4)$ acts on $H^3(T,\integers)$ as an
anti-involution with respect to the $\eisenstein$-module structure.  Let
$\chi$ be the corresponding integral anti-involution $i\comp \kappa^* \comp i^{-1}$
of $\chyper^4$.  This associates to $S$ and a choice of $i$
a pair $(L,\chi)\in K_0$, and defines a period map
\begin{equation}
  \moduli_0^\reals \map P\Gamma\backslash K_0\;,
\end{equation}
which we show is an isomorphism of real analytic orbifolds.

Another way to look at $K_0$ is as a disjoint union of incomplete real
hyperbolic manifolds.  To see this, let $\rhyper^4_\chi$ be
the set of fixed points in $\chyper^4$ of $\chi$.  Then
\[
    K_0 = \coprod_{\chi} ( \rhyper^4_\chi - {\calH} )\;,
\]
where $\chi$ varies over the integral anti-involutions of $\chyper^4$.
Now let $C$ be a set of of representatives for the
conjugacy classes of integral anti-involutions of $\chyper^4$ under the
action of $P\Gamma$.
Let $P\Gamma_\chi$ be the centralizer of $\chi$
in $P\Gamma$.  Then the quotient of $K_0$ by
$P\Gamma$ is 
\[
P\Gamma\backslash K_0=
    \coprod_{\chi \in C} P\Gamma_\chi\backslash( \rhyper^4_\chi - {\calH} )\;.
\]
To understand this quotient in  detail, we need
to  classify the integral anti-involutions $\chi$ of $\chyper^4$, modulo
the action of $P\Gamma$.
One shows that there are just five classes,
given by
\begin{equation}
\chi_j(z_0, \ldots, z_4) 
  = (\bar z_0, \epsilon_1\bar z_1, \epsilon_2\bar z_2, \epsilon_3\bar z_3, \epsilon_4\bar z_4)\;,
\end{equation}
where $j$ of the $\epsilon_i$ are $-1$ and the rest are $+1$.  
It is clear that each $P\Gamma_{\chi_j}$ is a subgroup of the
projective automorphism group of the $\integers$-lattice
$\Lambda^{\chi_j}$ fixed by $\chi_j$, and one can check that it is the
full projective isometry group.  Computing the quadratic forms on the
$\Lambda^{\chi_j}$ leads to the quadratic forms used to describe the
$P\Gamma_j$ in \eqref{mainequation}, so $P\Gamma_{\chi_j}=P\Gamma_j$.  This
yields \eqref{mainequation}, where $\Delta_j=\rhyper^4_{\chi_j}\cap\calH$.  
We found the Coxeter
diagrams by using Vinberg's algorithm \cite{Vinberg}.  

In order to carry out the gluing process leading to Theorem~\ref{maintheorem}, we
computed which points of the Weyl chambers $C_j$ lie in $\calH$; it
turns out that $C_j\cap H$ is a union of faces of $C_j$.  Then we had
to figure out which faces of the $C_j$ and $C_{j'}$ to glue to each
other and how; for this we studied how the various $\rhyper^4_\chi$
meet in $\chyper^4$.  Finally we worked out the result of
the gluing by  explicitly manipulating polyhedra in $\rhyper^4$.


%
\end{document}